\documentclass[leqno,11pt]{amsart}
 
\usepackage{amsmath,amsthm}
\usepackage{amssymb}
\usepackage{euscript}
\usepackage[frame,ps,matrix,arrow,curve,rotate]{xy} 

\numberwithin{equation}{subsection}  

\newcommand{\sqsp}{\renewcommand{\baselinestretch}{1.1}\tiny\normalsize}

\newtheorem{thm}[subsection]{Theorem}
\newtheorem{prop}[subsection]{Proposition}
\newtheorem{cor}[subsection]{Corollary}
\newtheorem{lemma}[subsection]{Lemma}

\newtheorem*{thm-a}{Main Theorem}
\newtheorem*{cor-a}{Corollary}

\theoremstyle{definition}

\DeclareMathOperator{\End}{End}  
\DeclareMathOperator{\Hom}{Hom}
\DeclareMathOperator{\Ext}{Ext}
\DeclareMathOperator{\Obs}{Obs}

\begin{document}
\title{Deformation Theory of Modules}
\author{Donald Yau}

\begin{abstract}
Algebraic deformations of modules over a ring are considered.  The resulting theory closely resembles Gerstenhaber's deformation theory of associative algebras.
\end{abstract}

\email{dyau@math.uiuc.edu}
\address{Department of Mathematics, University of Illinois at Urbana-Champaign, Urbana, Illinois 61801, USA}


\maketitle

\sqsp

\section{Introduction}

In a series of papers, Gerstenhaber \cite{ger1,ger2,ger3,ger4} developed an algebraic deformation theory of associative algebras.  The idea is that given an associative algebra $R$, one deforms its multiplication while keeping the associative property and its underlying additive structure.  The infinitesimal of a deformation can then be identified with a $2$-cocycle in the Hochschild cohomology of $R$ with coefficients in itself.  As a result, the vanishing of this cohomology group implies that $R$ is rigid, in the sense that every deformation is equivalent to the original multiplication.  Related issues such as identifying obstructions to extending a $2$-cocycle to a deformation were also considered.  This theory has been extended in various directions.  For example, in \cite{ger5} Gerstenhaber and Wilkerson studied deformations of differential graded modules and algebras, in which only the differential is deformed.  The article \cite{fox} by Fox is a very readable introduction to the theory of algebraic deformations.

The purpose of this paper is to develop a deformation theory of modules over an associative ring.  Let $M$ be a (left) $R$-module, where $R$ is an associative, not-necessarily commutative $k$-algebra with $k$ a commutative ring with unit.  In section \ref{sec:infinitesimal} below we will show that the infinitesimal of a deformation of $M$ is a $1$-cocycle in the Hochschild cohomology of $R$ with coefficients in the $R$-$R$-bimodule $\End(M)$ of $k$-linear endomorphisms of $M$.  As a consequence of this, it is observed that the vanishing of the first Hochschild cohomology group of $R$ with coefficients in $\End(M)$ implies that $M$ is rigid (Corollary \ref{cor:rigid}).  In section \ref{sec:approx}, the obstructions to extending a $1$-cocycle to a deformation of $M$ are identified.  In Gerstenhaber's terminology, a $1$-cocycle that is extendible to a deformation is called \emph{integrable}.  We approach the integrability question by considering the slightly more general question of extending an \emph{approximate deformation of order} $n$ to one of higher order.  The obstructions are all $2$-cocycles (Lemma \ref{lem:obs}) and the vanishing of their cohomology classes is equivalent to extendibility of a $1$-cocycle to a deformation (Theorem \ref{thm:obs}).  The question of when two extensions are equivalent is also considered (Proposition \ref{prop:obs}).  Notice the similarity of our results with Gerstenhaber's, except for a dimension shift.

We should point out that algebraic deformations of modules were first considered by Donald and Flanigan \cite{df}.  Their theory has the advantage that the definition of a deformation is very similar to Gerstenhaber's.  On the other, they have to restrict themselves to finite dimensional algebras $R$ over a field $k$ and finite dimensional $R$-modules $M$.  Our theory has the (very) slight advantage that we do not need such assumptions on $k$, $R$, and $M$.

Moreover, our deformation of an $R$-module $M$ can be thought of as a deformation of the algebra morphism $R \to \End(M)$.  In particular, our theory is a special case of the deformation of an algebra morphism considered by Gerstenhaber and Schack \cite{gs1,gs2,gs4,schack}.  In those papers, they studied cohomology and deformations of presheaves of algebras, with an algebra map $B \to A$ being a special case.  Furthermore, our results are quite similar to the two main examples discussed by Gerstenhaber and Schack in  \cite{gs3}.  The first of these examples is due to Nijenhuis \cite{nij}, which concerns deformations of a subalgebra $B$ inside a fixed algebra $A$.  The second example is about deformations of a group representation $kG \to \End(V)$, in which the group $G$, and hence the algebra $kG$, is fixed.  Although our results can be regarded as special cases of those of Gerstenhaber and Schack, it seems worthwhile to record the results in details in this special case for future reference.

We present our theory entirely in terms of left modules, but one can adapt it easily to the cases of bimodules and right modules.  For instance, it is a standard fact that an $R$-$R$-bimodule $M$ is equivalent to a left $R \otimes_k R^{\text{op}}$-module, via the formula $(r \otimes s)m = rms$.

In the following section, we will first recall some basic definitions about Hochschild cohomology.  The reader who is familiar with Hochschild cohomology can safely skip this section and go directly to section \ref{sec:infinitesimal}.


\section{Hochschild cohomology}
\label{sec:Hoch}

The purpose of this section is to give a brief account of Hochschild cohomology that is relevant to this paper.  Much more detailed discussions can be found in many references, for example, Weibel \cite{weibel}.

Throughout this paper, we work with a fixed commutative ground ring $k$ with unit, an associative, not-necessarily commutative $k$-algebra $R$, and a (left) $R$-module $M$.  For convenience, we write $\otimes$ for $\otimes_k$, tensoring over $k$, and $R^{\otimes n}$ for $R \otimes \cdots \otimes R$ ($n$ factors).

The set of $k$-linear endomorphisms of $M$, written $\End(M)$, has a natural structure as an $R$-$R$-bimodule.  In fact, if $r$ and $s$ are elements in $R$, $g \in \End(M)$, and $m \in M$, then 
   \[
   (rgs)(m) ~=~ rg(sm).
   \]
Moreover, $\End(M)$ is also an associative $k$-algebra with composition of endomorphisms as product.

Consider the cochain complex $\Hom_k(R^{\otimes *}, \End(M))$ in which the $n$th dimension consists of the $k$-linear maps from $R^{\otimes n}$ to $\End(M)$.  The differential
   \[
   d_n \colon \Hom_k(R^{\otimes n}, \End(M)) \to \Hom_k(R^{\otimes n+1}, \End(M))
   \]
is given by the formula
   \begin{multline}
   \label{eq:d}
   d_nf(a_0 \otimes \cdots \otimes a_n) \\
     ~=~ a_0 f(a_1 \otimes \cdots \otimes a_n) 
     ~+~ \sum_{i=1}^n (-1)^i f(a_0 \otimes \cdots \otimes a_{i-1}a_i \otimes \cdots \otimes a_n) \\ 
     ~+~ (-1)^{n+1} f(a_0 \otimes \cdots \otimes a_{n-1})a_n.
   \end{multline}
It is a standard fact that $d^2 = 0$.  The cohomology groups of this cochain complex are the \emph{Hochschild cohomology groups of} $R$ \emph{with coefficients in} $\End(M)$, denoted by $H^*_k(R,\End(M))$.


\section{Deformations and Infinitesimals}
\label{sec:infinitesimal}

We remind the reader that we are working with a fixed ground ring $k$, an associative, not-necessarily commutative $k$-algebra $R$, and a (left) $R$-module $M$.

In this section, we will define a deformation and a formal automorphism of $M$.  Then we observe that the infinitesimal deformation is a $1$-cocycle in the Hochschild cohomology of $R$ with coefficients in $\End(M)$.  A consequence of this is that the vanishing of the first Hochschild cohomology group of $R$ with coefficients in $\End(M)$ implies that $M$ is rigid.

To define a deformation of $M$, notice that the $R$-module structure on $M$ is equivalent to a map
   \[
   \xi \colon R ~\to~ \End(M)
   \]
of associative $k$-algebras, where $\End(M)$ denotes the set of $k$-linear endomorphisms of $M$.  Indeed, the map $\xi$ is given by
   \[
   \xi(r)(m) ~=~ rm
   \]
for elements $r \in R$ and $m \in M$.

With this in mind, we now define a (\emph{formal}) \emph{deformation} of $M$ to be a formal power series
   \[
   \xi_t ~=~ \xi + t \xi_1 + t^2 \xi_2 + \cdots
   \]
in which each $\xi_i$ is a $k$-linear map from $R$ to $\End(M)$.  Moreover, it is required to be multiplicative, in the sense that
   \begin{equation}
   \label{eq:multi}
   \xi_t(rs) ~=~ \xi_t(r) \xi_t(s)
   \end{equation}
for all elements $r$ and $s$ in $R$.  The $k$-linear map $\xi_1$ is called the \emph{infinitesimal deformation} of $\xi_t$.

The multiplicative property is equivalent to the relations
   \begin{equation}
   \label{eq:relation}
   \xi_n(rs) ~=~ \sum_{i + j = n}\, \xi_i(r) \xi_j(s)
   \end{equation}
for $n \geq 0$, where we are using the convention that $\xi_0 = \xi$.  When $n = 0$, this is just restating that $\xi$ is multiplicative.  When $n = 1$, we obtain the relation
   \[
   \xi_1(rs) ~=~ \xi(r) \xi_1(s) + \xi_1(r) \xi(s),
   \]
or equivalently,
   \begin{equation}
   \label{eq:n=1}
   \xi(r)\xi_1(s) - \xi_1(rs) + \xi_1(r) \xi(s) ~=~ 0.
   \end{equation}
In other words, $\xi_1 \in \Hom_k(R, \End(M))$ is a $1$-cocycle (that is, in the kernel of $d_1$ \eqref{eq:d}) in the Hochschild cochain complex $\Hom_k(R^{\otimes *}, \End(M))$ of $R$ with coefficients in $\End(M)$.

Of course, one would like to regard the cohomology class of $\xi_1$ as the infinitesimal deformation.  In order to do this, we need an appropriate notion of equivalence of deformations.  Define a \emph{formal automorphism} of $M$ to be a formal power series
   \[
   \phi_t ~=~ 1 + t \phi_1 + t^2 \phi_2 + \cdots
   \]
in which each $\phi_i$ is a $k$-linear endomorphism of $M$ with $1$ denoting the identity map.  Two deformations $\xi_t$ and $\bar{\xi}_t$ of $M$ are said to be \emph{equivalent} if there exists a formal automorphism $\phi_t$ such that
   \begin{equation}
   \label{eq:equiv}
   \bar{\xi}_t ~=~ \phi_t^{-1} \xi_t \phi_t.
   \end{equation}
This clearly defines an equivalence relation on the set of deformations of $M$. 
We say that $M$ is \emph{rigid} if every deformation of $M$ is equivalent to $\xi$.

Two things should be pointed out about the equation \eqref{eq:equiv}.  First, the right-hand side of the equation is a composition of power series.  It is to be understood in the following sense.  Given $k$-linear endomorphisms $\phi$ and $\psi$ of $M$ and a $k$-linear map $g \colon R \to \End(M)$, 
   \[
   \psi g \phi \colon R ~\to~ \End(M)
   \]
is the map that sends an element $r \in R$ to the composite $k$-linear endomorphism $\psi g(r) \phi$ of $M$.  Second, if $\xi_t$ is a deformation of $M$ and $\phi_t$ is a formal automorphism, then the power series $\bar{\xi}_t$ defined by \eqref{eq:equiv} is a deformation of $M$.  In other words, the formal power series $\bar{\xi}_t$ so defined has the multiplicative property.  This is a trivial consequence of the multiplicative property of $\xi_t$ \eqref{eq:multi}.

With the notion of equivalence defined as above, the following result allows us to regard the infinitesimal deformation as a $1$-dimensional cohomology class which is well defined by the equivalence class of the deformation.

\medskip
\begin{prop}
\label{prop:well}
Given a deformation $\xi_t$ of $M$, the infinitesimal deformation $\xi_1$ is a $1$-cocycle in the Hochschild cochain complex of $R$ with coefficients in $\End(M)$.  Moreover, if $\bar{\xi}_t$ is another deformation of $M$ that is equivalent to $\xi_t$, then $\bar{\xi}_1 - \xi_1$ is a $1$-coboundary.
\end{prop}

\begin{proof}
The first assertion has already been established.  For the second assertion, we know that there exists a formal automorphism $\phi_t$ such that \eqref{eq:equiv} holds.  Expanding the right-hand side of that equation, we obtain
   \[
   \begin{split}
   \phi_t^{-1} \xi_t \phi_t 
     &~=~ (1 - t \phi_1 + \cdots \, )(\xi + t\xi_1 + \cdots \, )(1 + t \phi_1 + \cdots \, ) \\
     &~=~ \xi + t(\xi_1 + \lbrack \xi, \phi_1 \rbrack) + \text{ higher terms in } t.
   \end{split}
   \]
Here $\lbrack \xi, \phi_1 \rbrack$ denotes $\xi \phi_1 - \phi_1 \xi$, which is  a $1$-coboundary in the Hochschild cochain complex of $R$ with coefficients in $\End(M)$.  This proves the second assertion.
\end{proof}

Now we would like to have a cohomological criterion for the rigidity of $M$.

Notice that if $\xi_t$ is a deformation of $M$ of the form
   \[
   \xi_t ~=~ \xi + t^l \xi_l + t^{l+1} \xi_{l+1} + \cdots
   \]
for some $l \geq 1$ (that is, $\xi_1 = \cdots = \xi_{l-1} = 0$), then $\xi_l$ is a $1$-cocycle.  This is again derived from \eqref{eq:relation}.  Suppose in addition that $\xi_l$ is a $1$-coboundary so that $\xi_l = \lbrack \xi, \phi_l \rbrack$ for some $\phi_l \in \End(M)$.  Let $\phi_t$ be the formal automorphism
   \[
   \phi_t ~=~ 1 - t^l \phi_l.
   \]
Then $\xi_t$ is equivalent to the following deformation of $M$:
   \[
   \begin{split}
   \phi_t^{-1} \xi_t \phi_t
     &~=~ (1 + t^l \phi_l + \cdots \, )(\xi + t^l \xi_l + \cdots \, )(1 - t^l \phi_l) \\
     &~\equiv~ \xi + t^l(\xi_l - \lbrack \xi, \phi_l \rbrack) \pmod{t^{l+1}} \\
     &~\equiv~ \xi \pmod{t^{l+1}}.
   \end{split}
   \]
In other words, if the first nontrivial $\xi_l$ in the deformation $\xi_t$ is cohomologous to $0$, then $\xi_t$ is equivalent to a deformation in which the first nontrivial term (after $\xi$) is $t^{l+1}$.  We record this formally in the following result.

\medskip
\begin{thm}
\label{thm:coh}
Let $\xi_t$ be a deformation of $M$.  Then there exists an $l \in \lbrace 1, 2, \ldots, \infty \rbrace$ such that $\xi_t$ is equivalent to a deformation of the form
   \[
   \bar{\xi}_t ~=~ \xi + t^l \bar{\xi}_l + t^{l+1} \bar{\xi}_{l+1} + \cdots
   \]
in which $\bar{\xi}_l$ is not cohomologous to $0$.
\end{thm}

An immediate consequence of this result is a cohomological criterion for the rigidity of $M$.

\medskip
\begin{cor}
\label{cor:rigid}
If the Hochschild cohomology group $H^1_k(R,\End(M))$ is trivial, then $M$ is rigid.
\end{cor}

Here is an example of a rigid module.  Let $k$ be a field, $R$ be a finite dimensional, separable $k$-algebra, and $M$ be any (left) $R$-module.  Then there is an isomorphism
   \[
   H^i_k(R, \End(M)) ~\cong~ \Ext^i_{R/k}(M,M),
   \]
which is $0$ if $i > 0$.  (See Example 8.7.6, Lemma 9.1.9, and Theorem 9.2.11 in \cite{weibel} for a proof.)   In particular, $M$ is rigid.  Here $\Ext^i_{R/k}$ is the ``relative'' $\Ext$ group defined as follows.  Let $U$ be the forgetful functor from $R$-modules to $k$-modules and $F$ be the functor $- \otimes_k R$ going the other way.  Then $F$ is a left adjoint to $U$, and therefore $FU$ is a cotriple on the category of $R$-modules.  Given an $R$-module $M$, we obtain in the usual way a simplicial $R$-module $(FU)^*M$ by applying the cotriple $FU$ repeatedly.  The relative $\Ext$ groups are then the cohomology groups of the cochain complex associated to the cosimplicial $R$-module $\Hom_R((FU)^*M, M)$.


\section{Extending approximate deformations}
\label{sec:approx}

We continue to use the notations from the previous section.  In particular, $R$ is an associative, not-necessarily commutative algebra over a commutative ground ring $k$ and $M$ is a (left) $R$-module.  The $R$-module structure on $M$ is considered as a map $\xi \colon R \to \End(M)$ of associative $k$-algebras.

In the previous section, we observed that the infinitesimal deformation is a $1$-cocycle in the Hochschild cochain complex of $R$ with coefficients in $\End(M)$.  It is not true that an arbitrary $1$-cocycle $\sigma$ is the infinitesimal deformation of a deformation of $M$.  Following Gerstenhaber \cite{ger1}, we say that a $1$-cocycle $\sigma$ is \emph{integrable} if it is the infinitesimal deformation of a deformation of $M$.  The purpose of this section is to identify the obstructions to the integrability of a $1$-cocycle.

We will actually consider a slightly more general problem.  Following Gerstenhaber and Wilkerson \cite{ger5}, define an \emph{approximate deformation of order} $m$ ($\geq 1$) to be a formal power series 
   \[
   \xi_t ~=~ \xi + t \xi_1 + t^2 \xi_2 + \cdots + t^m \xi_m
   \]
in which each $\xi_i$ is a $k$-linear map from $R$ to $\End(M)$.  It is required to satisfy \eqref{eq:multi} $\pmod{t^{m+1}}$.  Equivalently, it satisfies \eqref{eq:relation} for $n \leq m$.  A formal automorphism is defined just as before.  If $\bar{\xi}_t$ is another approximate deformation of order $m$, then we say that $\xi_t$ and $\bar{\xi}_t$ are equivalent if there exists a formal automorphism $\phi_t$ such that 
   \[
   \bar{\xi}_t ~\equiv~ \phi_t^{-1} \xi_t \phi_t \pmod{t^{m+1}}.
   \]

Given a $1$-cocycle $\sigma$ in the Hochschild cochain complex of $R$ with coefficients in $\End(M)$, integrating $\sigma$ is equivalent to extending $\xi_t = \xi + t \sigma$, an approximate deformation of order $1$, to approximate deformations of all higher orders.  To this end, we now consider the problem of identifying the obstructions to extending an approximate deformation of order $m$ to one of order $m + 1$.

Let, then, $\xi_t = \sum_{i=0}^m\, t^i \xi_i$ be an approximate deformation of order $m$.  Consider the $k$-linear map $\Obs(\xi_t) \in \Hom_k(R \otimes R, \End(M))$ defined by
   \[
   \Obs(\xi_t)(a \otimes b) 
     ~\buildrel \text{def} \over =~ \sum_{i=1}^m\, \xi_i(a) \xi_{m+1-i}(b)
   \]
for all $a \otimes b \in R \otimes R$.  We call $\Obs(\xi_t)$ an \emph{obstruction cocycle}.  The terminology will be justified by the following two results.

\medskip
\begin{lemma}
\label{lem:obs}
The element $\Obs(\xi_t)$ is a $2$-cocycle in the Hochschild cochain complex of $R$ with coefficients in $\End(M)$.
\end{lemma}

\begin{proof}
In Hochschild cohomology, the differential $d_2$ sends an element $\varphi \in \Hom_k(R^{\otimes 2}, \End(M))$ to $d_2 \varphi \in \Hom_k(R^{\otimes 3}, \End(M))$ given by
   \[
   (d_2 \varphi)(a \otimes b \otimes c) 
    ~=~ \xi(a) \varphi(b \otimes c) - \varphi(ab \otimes c) + \varphi(a \otimes bc) - \varphi(a \otimes b) \xi(c)
   \]
for all $a \otimes b \otimes c \in R^{\otimes 3}$.

We perform the following computation:
   \[
   \begin{split}
   &d_2 \Obs(\xi_t)(a \otimes b \otimes c) \\
   &=~ \xi(a) \sum_{i=1}^m\, \xi_i(b) \xi_{m+1-i}(c) 
      ~-~ \sum_{i=1}^m\, \xi_i(ab) \xi_{m+1-i}(c) \\
   &\quad +~ \sum_{i=1}^m \, \xi_i(a)\xi_{m+1-i}(bc) 
      ~-~ \sum_{i=1}^m\, \xi_i(a)\xi_{m+1-i}(b)\xi(c) \\
   &=~ \xi(a) \sum_{i=1}^m\, \xi_i(b) \xi_{m+1-i}(c)
      ~-~ \sum_{i=1}^m\, \left(\sum_{j=0}^i \xi_j(a)\xi_{i-j}(b)\right)\xi_{m+1-i}(c) \\ 
   &\quad +~ \sum_{i=1}^m\, \xi_i(a)\left(\sum_{j=0}^{m+1-i}\, \xi_j(b)\xi_{m+1-i-j}(c)\right) 
      ~-~ \sum_{i=1}^m\, \xi_i(a)\xi_{m+1-i}(b)\xi(c).
   \end{split}
   \]
In this last step, call the sum of the first two terms $A$ and the sum of the last two terms $B$.  Then we have
   \[
   A ~=~ - \sum_{i=1}^m\, \left(\sum_{j=1}^i\, \xi_j(a)\xi_{i-j}(b) \right)\xi_{m+1-i}(c) 
     ~=~ - \sum_{\substack{i+j+l ~=~ m+1 \\i,\, l > 0}}\, \xi_i(a)\xi_j(b)\xi_l(c).
   \]
Similarly, we have
   \[
   B ~=~ \sum_{i=1}^m\, \xi_i(a) \left( \sum_{j=0}^{m-i}\, \xi_j(b)\xi_{m+1-i-j}(c)\right)
     ~=~ \sum_{\substack{i+j+l ~=~ m+1 \\i,\, l > 0}}\, \xi_i(a)\xi_j(b)\xi_l(c).
   \]
Therefore, $d_2 \Obs(\xi_t)(a \otimes b \otimes c) = A + B = 0$, as was to be shown.  

This finishes the proof of the lemma.
\end{proof}

With this in mind, we can now show that the cohomology class of the $2$-cocycle $\Obs(\xi_t)$ is exactly the obstruction to extending $\xi_t$ to an approximate deformation of one higher order.

\medskip
\begin{thm}
\label{thm:obs}
Let $\xi_t$ be an approximate deformation of order $m$.  Then $\xi_t$ extends to an approximate deformation of order $m + 1$ if, and only if, the $2$-cocycle $\Obs(\xi_t)$ is a $2$-coboundary.
\end{thm}

\begin{proof}
Suppose that an order $m+1$ extension $\bar{\xi}_t = \xi_t + t^{m+1} \xi_{m+1}$ exists.  Then for $a, b \in R$ we have that
   \[
   \begin{split}
   \xi_{m+1}(ab)
     &~=~ \sum_{i+j \,=\, m+1}\, \xi_i(a)\xi_j(b) \\
     &~=~ \xi_{m+1}(a)\xi(b) ~+~ \xi(a)\xi_{m+1}(b) ~+~ \sum_{i=1}^m\, \xi_i(a)\xi_{m+1-i}(b).
   \end{split}
   \]
Moving $\xi_{m+1}(ab)$ to the left-hand side, this is equivalent to the equation
   \[
   0 ~=~ (d_1 \xi_{m+1})(a \otimes b) + \Obs(\xi_t)(a \otimes b),
   \]
which implies that $\Obs(\xi_t)$ is a $2$-coboundary.

Since the argument above is reversible, this finishes the proof of the theorem.
\end{proof}

An immediate consequence of this result is an obstruction theoretic answer to the question of integrating a $1$-cocycle.

\medskip
\begin{cor}
\label{cor:int}
Let $\sigma \in \Hom_k(R, \End(M))$ be a $1$-cocycle.  Then there exists a sequence of $2$-cocycles, $\Obs_i$ $(i \geq 1)$, such that $\Obs_m$ is defined if and only if $\Obs_i$ for $i < m$ are all defined and are cohomologous to $0$.  Moreover, $\sigma$ is integrable if and only if $\Obs_i$ is cohomologous to $0$ for all $i \geq 1$.
\end{cor}

As a special case, we obtain a simple cohomological criterion which guarantees that an approximate deformation of any order can be extended to a deformation of $M$.

\medskip
\begin{cor}
\label{cor:ext}
If the Hochschild cohomology group $H^2_k(R,\End(M))$ is trivial, then every approximate deformation of order $m$ ($\geq 1$) extends to a deformation of $M$.
\end{cor}

Finally, we consider the question of whether an extension is unique up to equivalence.

Let $\xi_t = \xi + t \xi_1 + \cdots + t^m \xi_m$ be an approximate deformation of $M$ of order $m$ and let $\xi^\prime_t = \xi_t + t^{m+1}\xi^\prime_{m+1}$ and $\xi^{\prime\prime}_t = \xi_t + t^{m+1}\xi^{\prime\prime}_{m+1}$ be two order $m + 1$ extensions of $\xi_t$.  From the proof of Theorem \ref{thm:obs}, we know that
   \[
   d_1 \xi^\prime_{m+1} ~=~ - \Obs(\xi_t) ~=~  d_1 \xi^{\prime\prime}_{m+1}.
   \]
In particular, the difference $\xi^\prime_{m+1} - \xi^{\prime\prime}_{m+1}$ is a $1$-cocycle.

We now observe that the vanishing of the cohomology class of $\xi^\prime_{m+1} - \xi^{\prime\prime}_{m+1}$ is a sufficient condition for the two extensions to be equivalent.

\medskip
\begin{prop}
\label{prop:obs}
With the notations as above, the two extensions $\xi^\prime_t$ and $\xi^{\prime\prime}_t$ are equivalent if the $1$-cocycle $\xi^\prime_{m+1} - \xi^{\prime\prime}_{m+1}$ is a $1$-coboundary.
\end{prop}

\begin{proof}
The hypothesis says that there exists some $k$-linear endomorphism $\phi$ of $M$ such that
   \[
   \xi^{\prime\prime}_{m+1} ~=~ \xi^\prime_{m+1} ~+~ \lbrack \xi, \phi \rbrack.
   \]
Define a formal automorphism of $M$ by setting $\phi_t = 1 + t^{m+1}\phi$.  Then we have
   \[
   \begin{split}
   \phi^{-1}_t \xi^\prime_t \phi_t 
     &~=~ (1 - t^{m+1} \phi + \cdots \,)(\xi_t + t^{m+1}\xi^\prime_{m+1})(1 + t^{m+1}\phi) \\
   &~\equiv~ \xi_t ~+~ t^{m+1}(\xi^\prime_{m+1} ~+~ \lbrack \xi, \phi \rbrack) \pmod{t^{m+2}} \\
   &~\equiv~ \xi^{\prime\prime}_t \pmod{t^{m+2}}.
   \end{split}
   \] 
This shows that the two approximate deformations $\xi^\prime_t$ and $\xi^{\prime\prime}_t$ of order $m + 1$ are equivalent.
\end{proof}

\section{Acknowledgment}
The author thanks the referee for reading an earlier version of this paper and for pointing out the relationships between the results here and those in the literature, especially \cite{gs2,gs3,gs4,nij}.



\end{document}